\documentclass[11pt,a4paper]{article}
\usepackage[english]{babel}
\usepackage{amsmath}
\usepackage{amssymb}
\usepackage{amsfonts}
\usepackage{amsthm}
\usepackage{graphicx}
\usepackage{fancyhdr}

\newtheorem{theorem}{Theorem}

\newtheorem{proposition}[theorem]{Proposition}
\newtheorem{definition}[theorem]{Definition}
\newtheorem{conjecture}[theorem]{Conjecture}

\makeatletter
   \renewcommand{\section}{\@startsection{section}{1}{0mm}
   {\baselineskip}%
   {\baselineskip}{\normalfont\normalsize\scshape\centering}}%
   \makeatother

\fancypagestyle{plain}{%
  \fancyhf{}
  \fancyfoot[C]{\thepage}

}

\allowdisplaybreaks[1]
\begin{document}

\pagestyle{fancy}
\thispagestyle{plain}

\begin{center}
	\textsc{\Large On the Entropy of Random Fibonacci Words} \vspace{1.5ex}
		
	\textsc{Johan Nilsson}\vspace{1.5ex}
	
	{\small Centre for Mathematical Sciences, Lund University, LTH \\
		Box 118, SE-221 00 LUND, SWEDEN}
	
	{\small\texttt{johann@maths.lth.se}}
	
\end{center}

\begin{abstract}
The random Fibonacci chain is a generalisation of the classical Fibonacci substitution and is defined as the rule mapping $0\mapsto 1$ and $1 \mapsto 01$ with probability $p$ and $1 \mapsto 10$  with probability $1-p$ for $0<p<1$ and where the random rule is applied each time it acts on a $1$. We show that the topological entropy of this object is given by the growth rate of the set of inflated random Fibonacci words.
\end{abstract}

{\small
\noindent MSC2010 classification; 68R15 Combinatorics on words, 05A16 Asymptotic enumeration, 37B10 Symbolic dynamics.
}

\section{Introduction}

In \cite{godreche} Godr\`eche and Luck define the random Fibonacci chain by the generalised substitution 
\begin{equation}
\label{eq: def of Theta}
	\theta : 
	\left\{
	\begin{array}{rcl}
		0 & \mapsto & 1 \\
		1 & \mapsto &
		\left\{
		\begin{array}{rl}
			01 & \textnormal{with probability $p$} \\
			10 & \textnormal{with probability $1-p$}
		\end{array}
		\right.
	\end{array}
	\right.
\end{equation}
for $0<p<1$ and where the random rule is applied each time $\theta$ acts on a $1$. They introduce the random Fibonacci chain when studying quasi-crystalline structures and tilings in the plane. In their paper it is claimed without proof that the topological entropy of the random Fibonacci chain is given by the growth rate of the set of inflated random Fibonacci words. We give here a proof of this fact.

Before we can state our main theorem we need to introduce some notation. A word $w$ over an alphabet $\Sigma$ is a finite sequence $w_1w_2\ldots w_n$ of symbols from $\Sigma$. We let here $\Sigma = \{0,1\}$. We denote a sub-word or a factor of $w$ by $w[a,b] = w_a w_{a+1}w_{a+2}\ldots w_{b-1}w_b$ and similarly we let $W[a,b] =\{w[a,b]: w\in W\}$. By $|\cdot|$ we mean the length of a word and the cardinality of a set. Note that $|w[a,b]| = b-a+1$.

For two words $u = u_1u_2u_3\ldots u_n$ and $v = v_1v_2v_3\ldots v_m$ we denote by $uv$ the concatenation of the two words, that is, $uv = u_1u_2u_3\ldots u_n v_1 v_2 \ldots v_m$. Similarly we let for two sets of words $U$ and $V$ their product be the set $UV = \{uv: u\in U, v\in V\}$ containing all possible concatenations.

Letting $\theta$ act on the word $0$ repeatedly yields us an infinite sequence of words $r_n = \theta^{n-1}(0)$. We know that $r_1=0$ and $r_2=1$. But $r_3$ is either $01$ or $10$, with probability $p$ or $1-p$ respectively. The sequence $\{r_n\}_{n=1}^{\infty}$ converges in distribution to an infinite random word $r$. We say that $r_n$ is an inflated word (under $\theta$) in generation $n$ and we introduce here sets that correspond to all inflated words in generation $n$; 

\begin{definition}
\label{def: recursive def An}
Let $A_0 := \emptyset$, $A_{1} := \{ 0 \}$ and $A_{2} := \{ 1 \}$ and for $n\geq3$ we define recursively 
\begin{equation}
\label{eq: rec def An}
	A_n := A_{n-1}A_{n-2}\cup A_{n-2}A_{n-1}
\end{equation}
and we let $A := \lim_{n\to\infty}A_n$. 
\end{definition}

We shall later on, as a direct consequence of Proposition \ref{prop: equal prefix}, see that the set $A$ is a well defined set. It is clear from the definition of $A_n$ that all elements in $A_n$ have the same length, that is, for all $x,y \in A_n$ we have $|x| = |y|$. We shall frequently us the length of the elements in $A_n$, and therefore we introduce the notation 
\[
	f_n := |x| \quad \textnormal{for}\quad x\in A_n,\quad n\geq1
\]
and for completeness we set $f_0:=0$. From the recursion (\ref{eq: rec def An}) and the definition of $f_n$ we have immediately the following proposition 

\begin{proposition}
The numbers $f_n$ are the Fibonacci-numbers, that is, $f_0 =0$,  $f_1 = 1$ and $f_n = f_{n-1} + f_{n-2}$ for $n\geq 2$.
\end{proposition}

For a word $w$ we say that $x$ is a factor of $w$ if there are two words $u,v$ such that $w= uxv$. The factor set $F(S,n)$, of a set of words $S$, is the set of all factors of length $n$ of the words in $S$.  
We introduce the abbreviated notation $F_n$ for the special factor sets of the random Fibonacci words,
\[	
	F_n := F(A,f_n).
\]
It seems to be a hard problem to give an explicit expression for the size of $|F_n|$. We have to leave that question open, but we shall give a rough upper bound of $|F_n|$. 

The topological entropy of the random Fibonacci chain is defined as the limit $\lim_{n\to\infty} \frac{1}{n}\log |F(A,n)|$. The existence of this limit is direct by Fekete's lemma \cite{fekete} since we have sub-additivity,  $\log|F(A,n+m)|\leq\log|F(A,n)|+\log|F(A,m)|$.
We can now state the main result in this paper

\begin{theorem}
\label{thm: lim An = lim Fn}
The logarithm of the growth rate of the size of the set of inflated random Fibonacci words equals the topological entropy of the random Fibonacci chain, that is 
\begin{equation}
\label{eq: lim An = lim Fn}
	\lim_{n\to\infty}\frac{\log|A_n|}{f_n} 
	= \lim_{n\to\infty}\frac{\log|F_n|}{f_n}.
\end{equation}
\end{theorem}

The outline of the paper is that we start by studying the set $A_n$ and restate some of the results obtained by Godr\`eche and Luck. Thereafter we continue by looking at sets of prefixes of the set $A_n$. Next we give a finite method for finding the factor set $F_n$ and finally we present an estimate of $|F_n|$ in terms of $|A_n|$, leading to the proof of Theorem \ref{thm: lim An = lim Fn}.

\section{Inflated Words}

In this section we present the sets of inflated random Fibonacci words and give an insight to their structure. For some small values of $n$ we have directly from the definition 
\[	
A_1=\{0\},\quad  A_2  = \{1\}, \quad  A_3 = \{01,10\}, 
\quad A_4 = \{011,101, 110\},
\]
\[	A_5  = 
	\left\{ \begin{array}{c}
				01011, 01101, 01110, 10011, \\
				10101, 10110, 11001, 11010\phantom{,}
			\end{array}
	\right\}.
\]
For a word $x$ we denote the reversed word of $x$ by $x^r$, that is,  if
$x = x_1x_2\ldots x_n$ then $x^r = x_nx_{n-1}\ldots x_1$.
Similarly we write $A_n^r = \{x^r : x\in A_n\}$.

\begin{proposition}
For $n\geq1$ we have $A_n = A_n^r$.
\end{proposition}
\begin{proof}
We give a proof by induction on $n$. The basis cases $n = 1$ and $n = 2$ are clear as $A_1 = \{0\}$ and $A_2=\{1\}$. Now assume for induction that the statement of the proposition holds for $2\leq  n\leq p$. Using (\ref{eq: rec def An}) and the induction assumption give
\begin{align*}
	A_{p+1}^r 
		&= (A_pA_{p-1}\cup A_{p-1}A_{p})^r  \\
		&= A_{p-1}^rA_{p}^r\cup A_{p}^rA_{p-1}^r  \\
		&= A_{p-1}A_{p}\cup A_{p}A_{p-1}  \\
		&= A_{p+1},
\end{align*}
which completes the induction.
\end{proof}

Now let us turn to the question of the size of $A_n$. 
In the recursive definition (\ref{eq: rec def An}) of $A_n$  we have that the overlap is
\[	
	A_{n}A_{n-1} \cap A_{n-1}A_{n} = A_{n-2}A_{n-3}A_{n-2}.
\]
Knowing this we obtain the recursion for $n\geq 3$
\begin{equation}
\label{eq: long rec An}
	|A_n| = 2|A_{n-1}||A_{n-2}|- |A_{n-2}|^2|A_{n-3}|,
\end{equation}
with $|A_0|=0$, $|A_1|=1$ and $|A_2|=1$. The recursion (\ref{eq: long rec An}) was given by Godr\`eche and Luck in \cite{godreche}. They also gave the following proposition, which we for completeness restate here with its short proof.  

\begin{proposition}[\cite{godreche}]
With $|A_1|=1$ and $|A_2|=1$ we have for $n\geq 3$
\label{prop: short rec An}
\begin{equation}
\label{eq: short rec An}
	|A_n| = \frac{n-1}{n-2}|A_{n-1}||A_{n-2}|.
\end{equation}
\end{proposition}

\begin{proof}
Let $a_n = |A_n| /(|A_{n-1}||A_{n-2}|)$. Then we can rewrite the recursion (\ref{eq: long rec An}) into
\[	a_n = 2 -\frac{1}{a_{n-1}}
\]
with the initial condition $a_3 = 2$. We claim that $a_n = \frac{n-1}{n-2}$. The proof of the claim is a straight forward induction. The basis is clear, and if we assume the claim being true for $n=p$ we get for $n=p+1$
\[	a_{p+1} = 2 - \frac{1}{a_{p}} =  2 - \frac{1}{\frac{p-1}{p-2}} = \frac{p}{p-1},
\]
which completes the induction and the proof of the proposition.
\end{proof}

Without explicitly stating it, the following proposition was also given by Godr\`eche and Luck in \cite{godreche}. We present it here, filling in the details.

\begin{proposition}[\cite{godreche}]
\label{prop: An explicit}
For $n\geq 3$ we have 
\begin{equation}
	\label{eq: An explicit}
	|A_n| = (n-1) \prod_{i=2}^{n-1} (n-i)^{f_{i-2}}. 
\end{equation}
\end{proposition}

\begin{proof}
We give a proof by induction on $n$. For the basis steps, $n=3$ and $n=4$ we have $|A_3| = 2\cdot 1^0 = 2$ and $|A_4| = 3\cdot 2^0\cdot 1^1 = 3$. Now assume for induction that (\ref{eq: An explicit}) holds for $4 \leq n \leq p$. Then from (\ref{eq: short rec An}) and the induction assumption we have 
\begin{align*}
	|&A_{p+1}| = \\
	& = \frac{p}{p-1}|A_{p}||A_{p-1}| \\
	& = \frac{p}{p-1} 
\left( (p-1) \prod_{i=2}^{p-1} (p-i)^{f_{i-2}} \right)
\left( (p-2) \prod_{i=2}^{(p-1)-1} ((p-1)-i)^{f_{i-2}} \right) \\
	& = p \, (p-2) 
\left( \prod_{i=2}^{p-1} (p-i)^{f_{i-2}} \right)
\left( \prod_{i=2}^{(p-1)-1} (p-(i+1))^{f_{(i+1)-3}} \right) \\
	& = p \, (p-2) 
\left( \prod_{i=2}^{p-1} (p-i)^{f_{i-2}} \right)
\left( \prod_{i=3}^{p-1} (p-i)^{f_{i-3}} \right) \\
	& = p \, (p-2)^{f_1} (p-2)^{f_0}
\left( \prod_{i=3}^{p-1} (p-i)^{f_{i-2}+f_{i-3}} \right) \\
	& = p \, (p-1)^{f_0} (p-2)^{f_1}
\left( \prod_{i=3}^{p-1} (p-i)^{f_{i-1}} \right) \\
	& = p \, (p-1)^{f_0} (p-2)^{f_1}
\left( \prod_{i=4}^{p} (p-(i-1))^{f_{(i-1)-1}} \right) \\
	& = p \left( \prod_{i=2}^{(p+1)-1} ((p+1)-i)^{f_{i-2}} \right),
\end{align*}
which completes the induction.
\end{proof}

The sequence $\{|A_n|\}_{n=1}^{\infty}$ (see Table \ref{table: A-F-c} on page \pageref{table: A-F-c}) is the sequence A072042 in the On-line Encyclopedia of Integer Sequences \cite{oleis}.

From Proposition \ref{prop: An explicit} we have 
\begin{align}
	\frac{\log |A_{n}|}{f_{n}} 
		&=	\frac{1}{f_{n}}\log \left((n-1) \prod_{i=2}^{n-1} (n-i)^{f_{i-2}}\right) \nonumber\\
		&=	\frac{1}{f_{n}}\log (n-1) + \sum_{i=2}^{n-1} \frac{f_{i-2}}{f_n}\log(n-i).
	\label{eq: fib-log sum}
\end{align}
For large $n$ the Fibonacci number $f_n$ can be approximated by $\frac{1}{\sqrt{5}}(\frac{1+\sqrt{5}}{2})^n$. Hence the sum in (\ref{eq: fib-log sum}) converges when letting $n$ tend to infinity and therefore clearly implies the existence of the limit. A numerical calculation in Maple gives
\[	
	\lim_{n\to\infty} \frac{\log |A_{n}|}{f_{n}}   \approx  0.444399  \approx  \log 1.559553. 
\]

\section{Prefixes}

In this section we look at and present some properties of sets of prefixes of the inflated random Fibonacci words. 

\begin{proposition}
\label{prop: equal prefix}
For $n\geq3$ and $k\geq0$ we have 
\begin{equation}
	\label{eq: An prefix 1}	
	A_{n}[1,f_n-1] = A_{n+k}[1,f_n-1]
\end{equation}
and symmetrically 
\begin{equation}
\label{eq: An prefix 2}	
	A_{n}[2,f_n] = A_{n+k}[f_{n+k}-f_{n}+2,f_{n+k}].
\end{equation}
\end{proposition}

\begin{proof}
We put our initial attention to equality (\ref{eq: An prefix 1}). Let us first consider the case with $k=1$. We will give a proof by induction on $n$. For the basis we have 
\[
	A_3[1,1]= \{0,1\} = A_4[1,1] \quad\textnormal{and} \quad  A_4[1,2]= \{01,10,11\} = A_5[1,2].
\]
Now assume that (\ref{eq: An prefix 1}) holds for $4\leq n\leq m$. In the induction step $n=m+1$ we get
\begin{align*}
	A_{m+1}[1,f_{m+1}-1] 
		&\subset  (A_{m+1}A_{m} \cup A_{m}A_{m+1})[1,f_{m+1}-1] \\
		&= A_{m+2}[1,f_{m+1}-1].
\end{align*}
For the reversed inclusion let $x\in A_{m+2}[1,f_{m+1}-1]$. By the recursion identity (\ref{eq: rec def An}) we have to look at two cases. First, if $x\in (A_{m+1}A_m)[1,f_{m+1}-1]$ then we clearly have $x\in A_{m+1}[1,f_{m+1}-1]$ and we are done. 
Secondly, we have to deal with the case $x\in (A_{m}A_{m+1})[1,f_{m+1}-1]$. For this case we use the induction assumption twice to obtain
\begin{align*}
	A_{m-1}[1,f_{m-1}-1]
		&= A_{m}[1,f_{m-1}-1] \\
		&= \big(A_{m}[1,f_m-1]\big)[1,f_{m-1}-1] \\
		&= \big(A_{m+1}[1,f_m-1]\big)[1,f_{m-1}-1] \\
		&= A_{m+1}[1,f_{m-1}-1].
\end{align*}
By the help of the above equality it follows that
\begin{align*}
	\big(A_{m}A_{m+1}\big)[1,f_{m+1}-1] 
		&= A_m \big(A_{m+1}[1,f_{m-1}-1]\big) \\
		&=  A_m \big(A_{m-1}[1,f_{m-1}-1]\big) \\
		&= \big(A_{m}A_{m-1}\big)[1,f_{m+1}-1] \\
		&\subset \big(A_{m}A_{m-1}\cup A_{m-1}A_{m}\big)[1,f_{m+1}-1]\\
		&= A_{m+1}[1,f_{m+1}-1],
\end{align*}
which completes the induction. 

Now let us turn to the case with $k\geq 1$. We shall also prove this by induction. We proved the basis step, $k=1$, above. Now we assume that (\ref{eq: An prefix 1}) holds for $k=p$. For $k=p+1$ we have by what we proved above
\[	
	A_{n+p}[1,f_{n+p}-1] =  A_{n+p+1}[1,f_{n+p}-1].
\]
But then also $A_{n+p}[1,f_{n}-1] =  A_{n+p+1}[1,f_{n}-1]$ and therefore by our induction assumption 
\[	
	A_{n}[1,f_{n}-1] =  A_{n+p}[1,f_{n}-1] =  A_{n+p+1}[1,f_{n}-1],
\]
which completes the proof of (\ref{eq: An prefix 1}). Equality (\ref{eq: An prefix 2}) is direct from (\ref{eq: An prefix 1}) by symmetry. 
\end{proof}

\begin{proposition}
\label{prop: An sub An-1An-2}
For $n\geq4$ we have 
\begin{equation}
\label{eq: An sub An-1An-2}
	A_{n} \subset \big(A_{n-1}[1,f_{n-1}-1]\big) \{0,1\}^2 \big(A_{n-2}[2,f_{n-2}]\big).
\end{equation}
\end{proposition}

\begin{proof}
Denote the right hand side set in (\ref{eq: An sub An-1An-2}) by $B_n$. If $n=4$ we have 
$B_4 = \big(A_{3}[1,1]\big) \{0,1\}^2 = \{0,1\}^3$, which clearly is a superset of $A_4$. 
If $n\geq 5$ then by Proposition \ref{prop: equal prefix} we have $A_n[1,f_{n-1}] = B_n[1,f_{n-1}]$ and $A_n[f_{n-1}+2,f_n] = A_{n-2}[2,f_{n-2}] =  B_n[f_{n-1}+2,f_n]$. As $B_n$ by definition consist of all words starting with a prefixes in $A_n[1,f_{n-1}]$ and ending with a suffix in $A_n[f_{n-1}+2,f_n]$ it must contain all sequences from $A_n$.
\end{proof}

Note that there clearly is a symmetric analogous of the result of Proposition \ref{prop: An sub An-1An-2}.

\section{Factor Sets}

The aim of this section is to give a finite method for finding the factor set $F_n$. For some small values of $n$ we have 
\[	F(A_3,f_2) = \{0,1\}, \qquad	F(A_4,f_3) = \{01,10,11\}, 
\]
\[	F(A_5,f_4) = \{001,010,011,100,101,110,111\},
\]
\[	\begin{array}{rcl}
		F(A_6,f_5) &=& 
		\left\{ \begin{array}{l}
			00101,00110,00111,01001,01010,	\\
			01011,01100,01101,01110,01111,	\\
			10010,10011,10100,10101,10110,	\\
			10111,11001,11010,11011,11100,	\\
			11101,11110
				\end{array}
		\right\}.
	\end{array}
\]
We shall shortly see that for $n\geq 4$ the above sets $F(A_{n+1},f_n)$ actually is $F_n$. It is clear from above that $F_3 \neq F(A_4,f_3)$ as words in $F(A_5,f_4)$ have the sub-word $00$, while this sub-word is not in $F(A_4,f_3)$. By computer calculation we can find the size of $F(A_{n+1},f_n)$ for some additional small values of $n$, see Table \ref{table: A-F-c} on page \pageref{table: A-F-c}.

\begin{proposition}
\label{prop: FAn = FA(n+k)}
For $k\geq 1$ and $n\geq 4$ we have 
\begin{equation}
\label{eq: FAn = FA(n+k)}
	F(A_{n+1},f_n) = F(A_{n+k},f_n).
\end{equation}
\end{proposition}

\begin{proof}
We give a proof by induction on $k$. For the basis step, $k=2$, we see that $F(A_{4+1},f_4) = F(A_{4+2},f_4)$, (this is direct as the word $000$ is not a sub-word of an element in $F(A_{6},f_5)$). Hence we may assume that $n\geq 5$. From the recursion (\ref{eq: rec def An}) we have directly that $F(A_{n+1},f_n) \subset F(A_{n+2},f_n)$. For the reversed inclusion we have by the recursion identity (\ref{eq: rec def An}) and symmetry that it is enough to show that any sub-word of length $f_n$ in $A_{n+1}A_n$ is also a sub-word in $A_{n+1}$. 
Let us define $x_i = x[i,i-1+f_n]$ where $x\in A_{n+1}A_n$ and $1\leq i \leq f_{n+1}+1$. 

For $1\leq i \leq f_{n-1}+1$ we have directly that $x_i$ is a sub-word in $A_{n+1}$.

For $f_{n-1} +2 \leq  i  < 2 f_{n-1}$  we have by the definition of $x_i$ that it is a sub-word in $A_{n+1}A_{n}[f_{n-1}+2, f_{n+1}+f_{n-1}-1]$. (The inequalities defining the case are valid if $f_{n-1}\geq 3$, which requires that $n\geq 5$). By Proposition \ref{prop: equal prefix} and the recursion (\ref{eq: rec def An}) we have 
\begin{align*}
	\big(A_{n+1}A_{n}\big)[f_{n-1}+2,& \,f_{n+1}+f_{n-1}-1] = \\
		& = \big(A_{n+1}[f_{n-1}+2,f_{n+1}]\big)\big(A_{n}[1,f_{n-1}-1]\big) \\
		& = \big(A_{n}[2,f_{n}]\big)\big(A_{n-1}[1,f_{n-1}-1]\big) \\
		& = \big(A_{n}A_{n-1}\big)[2,f_{n+1}-1] \\
		& \subset \big(A_{n}A_{n-1}\cup A_{n-1}A_{n}\big)[2,f_{n+1}-1] \\
		& = A_{n+1}[2,f_{n+1}-1],
\end{align*}
which gives that $x_i$ is a sub-word of $A_{n+1}$.

For $2 f_{n-1} \leq  i  \leq f_{n+1}$ we have that $x_i$ is a sub-word in $\big(A_{n+1}A_{n}\big)[f_{n}+2, f_{n+2}]$. Proposition \ref{prop: equal prefix} and the recursion (\ref{eq: rec def An}) now gives 
\begin{align*}
	\big(A_{n+1}A_{n}\big)[f_{n}+2,f_{n+2}] 
		& =   \big(A_{n+1}[f_{n}+2,f_{n+1}]\big) A_{n} \\
		& =   \big(A_{n-1}[2,f_{n-1}]\big)A_{n} \\
		& =   \big(A_{n-1}A_{n}\big)[2,f_{n+1}-1] \\
		& \subset  \big(A_{n}A_{n-1}\cup A_{n-1}A_{n}\big)[2,f_{n+1}-1] \\
		& =  A_{n+1}[2,f_{n+1}-1],
\end{align*}
and again we see that $x_i$ is a sub-word of $A_{n+1}$.

For finally $i = f_{n+1}+1$ we have that $x_i$ is an element in $A_n$, and therefore also a sub-word in $A_{n+1}$, which concludes the proof of the basis step.

Assume for induction that equality (\ref{eq: FAn = FA(n+k)}) holds for $k=p$. For the induction step, $k=p+1$, we have from what we just proved
\[
	F(A_{n+p},f_{n+p}) = F(A_{n+p+1},f_{n+p}).
\]
But then we must also have $F(A_{n+p},f_{n}) = F(A_{n+p+1},f_{n})$. The induction assumption now gives $F(A_{n},f_{n}) = F(A_{n+p},f_{n}) = F(A_{n+p+1},f_{n})$, which proves the proposition.
\end{proof}

\begin{proposition}
\label{prop: Fn = Fan}
	For $n\geq 4$ we have $F_n = F(A_{n+1},f_n)$.
\end{proposition}

\begin{proof}
It is clear that $F(A_{n+1},f_n) \subset F_n$. For the converse let $x\in F_n$. Then there are two words $u$ and $v$ such that $uxv\in A_m$ for some $m\geq n+1$, that is $x\in F(A_m,f_n)$. Proposition \ref{prop: FAn = FA(n+k)} now gives $F(A_m,f_n) = F(A_{n+1},f_n)$.
\end{proof}

\section{Upper Bound}

The aim of this section is to give an upper bound of the size of $|F_n|$ in terms of $|A_n|$. 

\begin{proposition}
\label{prop: An][An < cAn}
For $n\geq 3$ and $1\leq k\leq f_n-1$ we have 
\begin{equation}
\label{eq: An][An < cAn}	
	\big|A_n[1,k]\big|\cdot\big|A_n[k+1,f_n]\big| \leq 4^{n-2}|A_n|.
\end{equation}
\end{proposition}

\begin{proof}
It is clear from the symmetric structure of $A_n$ that we only have to give a proof for the case $1\leq k\leq \lfloor\frac{1}{2}f_n\rfloor$. We give a proof by induction on $n$. For the basis case $n=3$ we have
\[	
	|A_3[1,1]|\cdot|A_3[2,2]| = |\{0,1\}|\cdot|\{0,1\}| = 2\cdot 2 \leq 4\cdot 2  = 4^{3-2}|A_3|.
\]
Similarly we have for $n=4$
\[	
	|A_4[1,1]|\cdot|A_4[2,3]| \leq  |A_4|\cdot |A_4| = 3 \cdot |A_4|  \leq 4^{4-2}|A_4|
\]
and for $n=5$ with $k=1, 2$
\[
	|A_5[k,1]|\cdot|A_5[k+1,5]| \leq  |A_5|\cdot |A_5| = 8 \cdot |A_5|  \leq 4^{5-2}|A_5|.
\]
Now assume for induction that (\ref{eq: An][An < cAn}) holds for $5\leq n\leq p$. For the induction step, $n=p+1$, we have from Proposition \ref{prop: equal prefix}, Proposition \ref{prop: An sub An-1An-2}, the induction assumption and (\ref{eq: short rec An}) that
\begin{align*}
	\big|A_{p+1}&[1,k]\big| \cdot \big|A_{p+1}[k+1,f_{p+1}]\big| = \\
		& = \big|A_{p}[1,k]\big| \cdot \big|A_{p+1}[k+1,f_{p+1}]\big| \\
		& \leq \big|A_{p}[1,k]\big| \cdot \left| \Big(\big(A_{p}[1,f_p-1]\big) \{0,1\}^2 \big(A_{p-1}[2,f_{p-1}]\big)\Big)[k+1,f_{p+1}]\right| \\
		& = \big|A_{p}[1,k]\big| \cdot \big|\big(A_{p}[k+1,f_p-1]\big) \{0,1\}^2 \big(A_{p-1}[2,f_{p-1}]\big)\big| \\
		& \leq \big|A_{p}[1,k]\big| \cdot \big|A_{p}[k+1,f_{p}]\big| \cdot 2^2 \cdot \big|A_{p-1}[2,f_{p-1}]\big| \\
		&\leq 4^{p-2} \cdot |A_p|\cdot 4 \cdot |A_{p-1}| \\
		&\leq 4^{(p+1)-2} |A_{p+1}|,
\end{align*}
which completes the induction.
\end{proof}

The estimate in Proposition \ref{prop: An][An < cAn} seems to be far off from the true upper bound of the left hand side in (\ref{eq: An][An < cAn}). 
If we let 
\[	
	c_n = \max_{1\leq k \leq f_n-1} \frac{\big|A_n[1,k]\big|\cdot\big|A_n[k+1,f_n]\big|}{|A_n|}
\]
then computer calculations (see Table \ref{table: A-F-c} on page \pageref{table: A-F-c}) motivate the following conjecture

\begin{conjecture}
\label{conj: An][An < cAn}
There is a constant $C$ such that for all $n\geq 3$ and $1\leq k\leq f_n-1$ we have 
\[	
	\big|A_n[1,k]\big|\cdot\big|A_n[k+1,f_n]\big| \leq C \cdot |A_n|.
\]
\end{conjecture}
From our computer calculations (see Table \ref{table: A-F-c} on page \pageref{table: A-F-c}) we further conjecture that we may find $C<3$, but we have to leave this question open as well as the more delicate question of finding the optimal value of such a constant, $ c_o = \sup_{n\geq 3} c_n$.

We can now give the estimate we set out to find. Compared to computer calculations our estimate seems to be far from the optimal one, but it will be sufficient for our purpose.

\begin{proposition}
\label{prop: F_n < A_n}
For $n\geq 3$ we have 
\[	
	|F_n| \leq 2 \left(4^{n-2} f_{n-1}+1\right) |A_n|.
\]
\end{proposition}

\begin{proof}
For $n=3$ we have 
\[	
	|F_3| = 4 \leq 2 (4 \cdot 1 +1)\cdot 2 = 2(4^{3-2} f_{3-1}+1)|A_3|.
\]
For $n \geq 4 $ we have from Proposition \ref{prop: Fn = Fan} that $F_n = F(A_{n+1},f_n)$. Therefore it is enough to give an estimarte of the number of subwords of length $f_n$ there are in the words in $A_{n+1}$. Let us start by considering the number of subwords in the subset $A_{n}A_{n-1}$. We define $x_k = x[k,k-1+f_n]$, where $x\in A_{n}A_{n-1}$ for $1\leq k \leq f_{n-1}+1$. Clearly we have $x_k\in F_n$. If $k=1$ then we have $x_1 \in A_n$ and we see that there are $|A_n|$ different $x_1$'s. If $2\leq k\leq f_{n-1}+1$ we have by Proposition \ref{prop: equal prefix}
\[	
	x_k \in (A_n[k, f_n])(A_{n-1}[1, k-1]) = (A_n[k, f_n])(A_{n}[1, k-1]). 
\] 
Proposition \ref{prop: An][An < cAn} now gives that there are at most $4^{n-2} |A_n|$ different $x_k$'s for each of these $k$'s. Adding up the estimates and include the symmetric case, $x\in A_{n-1}A_{n}$, completes the proof.
\end{proof}

\section{Proof of the Theorem}

Finally, we have now gathered enough background to prove our theorem 

\begin{proof}[Proof of Theorem \ref{thm: lim An = lim Fn}]
We have already dealt with the existence of both limits in (\ref{eq: lim An = lim Fn}). Since clearly $A_n \subset F_n$ it now  follows from Proposition \ref{prop: F_n < A_n},
\begin{align*}
	\lim_{n\to\infty} \frac{\log|A_n|}{f_n} 
&\leq  \lim_{n\to\infty} \frac{\log|F_n|}{f_n} \\
&\leq  \lim_{n\to\infty} \frac{\log \big(2(4^{n-2} f_{n-1}+1) |A_n|\big)}{f_n} \\
\displaybreak[1]
&\leq  \lim_{n\to\infty} \frac{\log \big(4^{n-1}2^{n} |A_n|\big)}{f_n} \\
&\leq  \lim_{n\to\infty} \frac{ 3n \log 2}{f_n} +\lim_{n\to\infty} \frac{\log|A_n|}{f_n} \\
&=  \lim_{n\to\infty} \frac{\log |A_n|}{f_n},
\end{align*}
which completes the proof.
\end{proof}

\section{Numerics}

We present here the output from our computer calculations, which were obtained by a small JAVA-program on a standard PC.  
\newcounter{tablecounter}
\[
\refstepcounter{tablecounter}
\label{table: A-F-c}
\newlength{\mylen}
\settowidth{\mylen}{$|F(A_{n+1},f_n)|$}
\renewcommand{\arraystretch}{1.15}
\begin{tabular}{|r|r|r|r|r|r@{}l|}
	\multicolumn{7}{c}{\small Table \thetablecounter} \\
	\hline
	\multicolumn{1}{|c|}{$n$}&
	\multicolumn{1}{|c|}{$f_n$}&
	\multicolumn{1}{p{\mylen}|}{\centering $|A_n|$} & 
	\multicolumn{1}{p{\mylen}|}{\centering $|F_n|$}& 
	\multicolumn{1}{p{\mylen}|}{\centering $|F(A_{n+1},f_n)|$}  & 
	\multicolumn{2}{p{\mylen}|}{\centering $c_n$}  \\ \hline\hline
	\phantom{0}0& 0		& 0 		& 					& 			& \phantom{000}	&		\\ \hline
			1	& 1		& 1 		& 2					& 1			& 	&		\\ \hline 
			2	& 1		& 1 		& 2					& 2 		& 	&		\\ \hline
			3	& 2		& 2 		& 4					& 3 		& 2&.0		\\ \hline
			4	& 3		& 3 		& 7 				& 7 		& 2&.0		\\ \hline
			5	& 5		& 8 		& 22 				& 22		& 2&.0		\\ \hline
			6	& 8		& 30 		& 108 				& 108		& 2&.13333	\\ \hline
			7	& 13	& 288 		& 1356 				& 1356 		& 2&.11111	\\ \hline
			8	& 21	& 10080 	& 65800 			& 65800		& 2&.17143	\\ \hline
			9	& 34	& 3317760	& 30139200 			& 30139200	& 2&.16389	\\ \hline
\end{tabular}
\]

\section{Acknowledgment}

The author wishes to thank Michael Baake, Franz G\"ahler and Claudia L\"utke\-h\"olter
at Bielefeld University, Germany and Tomas Persson at the Polish Academy of Sciences, Warsaw, Poland for presenting the problem, our discussions of it and for reading drafts of the manuscript. This work was partly supported by the German Research Council (DFG), via CRC 701.

\end{document}